# Stochastic calculus for convoluted Lévy processes

CHRISTIAN BENDER[1] and TINA MARQUARDT[2]

[1]*Institute for Mathematical Stochastics, TU Braunschweig, Pockelsstr. 14, D-38106 Braunschweig, Germany. E-mail:* c.bender@tu-bs.de
[2]*Center of Mathematical Sciences, Munich University of Technology, D-85747 Garching, Germany. E-mail:* marquard@ma.tum.de

We develop a stochastic calculus for processes which are built by convoluting a pure jump, zero expectation Lévy process with a Volterra-type kernel. This class of processes contains, for example, fractional Lévy processes as studied by Marquardt [*Bernoulli* **12** (2006) 1090–1126.] The integral which we introduce is a Skorokhod integral. Nonetheless, we avoid the technicalities from Malliavin calculus and white noise analysis and give an elementary definition based on expectations under change of measure. As a main result, we derive an Itô formula which separates the different contributions from the memory due to the convolution and from the jumps.

## 1. Introduction

In recent years, fractional Brownian motion and other Gaussian processes obtained by convolution of an integral kernel with a Brownian motion have been widely studied as a noise source with memory effects (see, e.g., Alòs *et al.* [2], Bender [4], Biagini *et al.* [5] and the survey article by Nualart [15]). Potential applications for noise sources with memory are in such diverse fields as telecommunication, hydrology and finance, to mention a few.

In Marquardt [14], fractional Lévy processes were introduced. While capturing memory effects in a similar fashion as a fractional Brownian motion does, the convolution with a Lévy process provides more flexibility concerning the distribution of the noise (e.g., heavy tails). In this paper, we consider a larger class of processes by convolution of a rather general Volterra-type kernel with a centered pure jump Lévy process. These convoluted Lévy process may have jumps and/or memory effects depending on the choice of the kernel. Following the elementary $S$-transform approach developed by Bender [4] for fractional Brownian motion, we motivate and construct a stochastic integral with respect to convoluted Lévy processes. The integral is of Skorokhod type and so its zero expectation property makes it a possible choice to model an additive noise. As a main result, we derive an Itô formula for these integrals. The Itô formula clarifies the different influences of jumps and memory effects, which are captured in different terms.

The only other paper of which we are aware that treats integration for a similar class of processes is [7]. The class of filtered Poisson processes considered in their paper is analogously defined by replacing the Lévy process by a marked point process in the convolution. However, we emphasize that our approach allows the Lévy process to be of infinite variation and that our Itô formula for the Skorokhod integral is quite different from the one of Decreusefond and Savy [7] derive for a Stieltjes integral only.

The paper is organized as follows. After some preliminaries on Lévy processes and convoluted Lévy processes in Section 2, we discuss the $S$-transform in Section 3. The results from Section 3 motivate a definition for a Skorokhod integral with respect to convoluted Lévy processes which is given in Section 4. In this section, some basic properties of this integral are also discussed. Section 5 is devoted to the derivation of the Itô formula, while some results are specialized to fractional Lévy processes in Section 6.

# 2. Preliminaries

## 2.1. Basic facts on Lévy processes

We state some elementary properties of Lévy processes that will be needed below. For a more general treatment and proofs, we refer to Cont and Tankov [6] and Sato [18]. For notational convenience, we abbreviate $\mathbb{R}_0 = \mathbb{R} \setminus \{0\}$. Furthermore, $\|f\|$ is the ordinary $L^2$-norm of the function $f : \mathbb{R} \to \mathbb{R}$ and the corresponding inner product is denoted by $(f,g)_{L^2(\mathbb{R})}$. In this paper, we assume as given an underlying complete probability space $(\Omega, \mathcal{F}, P)$. Since the distribution of a Lévy processes $L$ on $(\Omega, \mathcal{F}, P)$ is infinitely divisible, $L$ is determined by its characteristic function in the Lévy–Khinchine form $\mathrm{E}[\mathrm{e}^{\mathrm{i}uL(t)}] = \exp\{t\psi(u)\}, t \geq 0$, where

$$\psi(u) = \mathrm{i}\gamma u - \tfrac{1}{2}u^2\sigma^2 + \int_{\mathbb{R}} (\mathrm{e}^{\mathrm{i}ux} - 1 - \mathrm{i}ux 1_{\{|x| \leq 1\}})\nu(\mathrm{d}x), \qquad u \in \mathbb{R}, \tag{1}$$

$\gamma \in \mathbb{R}$, $\sigma^2 \geq 0$ and $\nu$ is a Lévy measure on $\mathbb{R}$ that satisfies

$$\nu(\{0\}) = 0 \quad \text{and} \quad \int_{\mathbb{R}} (x^2 \wedge 1)\nu(\mathrm{d}x) < \infty.$$

For any measurable set $B \subset \mathbb{R}_0 \times \mathbb{R}$, let

$$N(B) = \sharp\{s \geq 0 : (L_s - L_{s-}, s) \in B\}$$

be the Poisson random measure on $\mathbb{R}_0 \times \mathbb{R}$, with intensity measure $n(\mathrm{d}x, \mathrm{d}s) = \nu(\mathrm{d}x)\,\mathrm{d}s$, that describes the jumps of $L$. Furthermore, let $\tilde{N}(\mathrm{d}x, \mathrm{d}s) = N(\mathrm{d}x, \mathrm{d}s) - \nu(\mathrm{d}x)\,\mathrm{d}s$ be the compensated jump measure of $L$ (see, e.g., Cont and Tankov [6], Definition 2.18).

Assuming that $\nu$ satisfies $\int_{|x|>1} x^2\nu(\mathrm{d}x) < \infty$, $L$ has finite mean and variance given by

$$\mathrm{var}(L(1)) = \int_{\mathbb{R}} x^2\nu(\mathrm{d}x) + \sigma^2. \tag{2}$$

If $\sigma = 0$ in (1), we call $L$ a *Lévy process without Brownian component*. In what follows, we will always assume that the Lévy process $L$ has no Brownian part. Furthermore, we suppose that $\mathrm{E}[L(1)] = 0$, hence $\gamma = -\int_{|x|>1} x\nu(\mathrm{d}x)$. Thus, (1) can be written in the form

$$\psi(u) = \int_{\mathbb{R}} (\mathrm{e}^{\mathrm{i}ux} - 1 - \mathrm{i}ux)\nu(\mathrm{d}x), \qquad u \in \mathbb{R}, \tag{3}$$

and $L = \{L(t)\}_{t\geq 0}$ can be represented as

$$L(t) = \int_0^t \int_{\mathbb{R}_0} x\tilde{N}(\mathrm{d}x, \mathrm{d}s), \qquad t \in \mathbb{R}. \tag{4}$$

In this case, $L$ is a martingale. In the sequel, we will work with a two-sided Lévy process $L = \{L(t)\}_{t\in\mathbb{R}}$, constructed by taking two independent copies $\{L_1(t)\}_{t\geq 0}$, $\{L_2(t)\}_{t\geq 0}$ of a one-sided Lévy process and setting

$$L(t) = \begin{cases} L_1(t), & \text{if } t \geq 0, \\ L_2(-t-), & \text{if } t < 0. \end{cases} \tag{5}$$

From now on, we will suppose that $\mathcal{F}$ is the completion of the $\sigma$-algebra generated by the two-sided Lévy process $L$ and will define $L^p(\Omega) := L^p(\Omega, \mathcal{F}, P)$.

## 2.2. Convoluted and fractional Lévy processes

We call a stochastic process $M = \{M(t)\}_{t\in\mathbb{R}}$ given by

$$M(t) = \int_{\mathbb{R}} f(t,s)L(\mathrm{d}s), \qquad t \in \mathbb{R}, \tag{6}$$

a *convoluted Lévy process* with kernel $f$. Here, $f : \mathbb{R} \times \mathbb{R} \to \mathbb{R}$ is a measurable function satisfying the following properties:

(i) $f(t,\cdot) \in L^2(\mathbb{R})$ for all $t \in \mathbb{R}$;
(ii) $f(t,s) = 0$ whenever $s > t \geq 0$, that is, the kernel is of Volterra type;
(iii) $f(0,s) = 0$ for almost all $s$, hence $M(0) = 0$.

Furthermore, we suppose that $L = \{L(t)\}_{t\in\mathbb{R}}$ is a Lévy process without Brownian component satisfying $\mathrm{E}[L(1)] = 0$ and $\mathrm{E}[|L(t)|^m] < \infty$ for all $m \in \mathbb{N}$. Hence, the process $M$ can be rewritten as

$$M(t) = \int_{\mathbb{R}} \int_{\mathbb{R}_0} f(t,s)x\tilde{N}(\mathrm{d}x, \mathrm{d}s), \qquad t \in \mathbb{R}. \tag{7}$$

Since $f(t,\cdot)\in L^2(\mathbb{R})$, the integral (7) exists in $L^2(\Omega,P)$ and

$$\mathrm{E}[M(t)^2]=\mathrm{E}[L(1)^2]\int_{\mathbb{R}} f^2(t,s)\,\mathrm{d}s=\mathrm{E}[L(1)^2]\|f(t,\cdot)\|^2_{L^2(\mathbb{R})}. \tag{8}$$

As an important class of examples of convoluted Lévy processes, we now consider univariate fractional Lévy processes. The name "fractional Lévy process" already suggests that it can be regarded as a generalization of fractional Brownian motion (FBM). We review the definition of a one-dimensional fractional Lévy process (FLP). For further details on FLP's, see Marquardt [13, 14].

***Definition 2.1 (Fractional Lévy Process (FLP)).*** *Let $L=\{L(t)\}_{t\in\mathbb{R}}$ be a Lévy process on $\mathbb{R}$ with $\mathrm{E}[L(1)]=0$, $\mathrm{E}[L(1)^2]<\infty$ and without Brownian component. For fractional integration parameter $0<d<0.5$, a stochastic process*

$$M_d(t)=\frac{1}{\Gamma(d+1)}\int_{-\infty}^{\infty}[(t-s)_+^d-(-s)_+^d]L(\mathrm{d}s),\qquad t\in\mathbb{R}, \tag{9}$$

*is called a* fractional Lévy process *(FLP).*

Note that the kernel (9) given by

$$f_t(s)=\frac{1}{\Gamma(1+d)}[(t-s)_+^d-(-s)_+^d],\qquad s\in\mathbb{R}, \tag{10}$$

satisfies conditions (i)–(iii). Thus, fractional Lévy processes are well defined and belong to $L^2(\Omega)$ for fixed $t$.

Moreover, the kernel can be represented by fractional integrals of the indicator function. Recall that for $0<\alpha<1$, the fractional integral of Riemann–Liouville type $I^\alpha_\pm$ is defined by

$$(I^\alpha_- f)(x)=\frac{1}{\Gamma(\alpha)}\int_x^\infty f(t)(t-x)^{\alpha-1}\,\mathrm{d}t,$$
$$(I^\alpha_+ f)(x)=\frac{1}{\Gamma(\alpha)}\int_{-\infty}^x f(t)(x-t)^{\alpha-1}\,\mathrm{d}t$$

if the integrals exist for almost all $x\in\mathbb{R}$.

In terms of these fractional operators, fractional Lévy processes can be rewritten as

$$M_d(t)=\int_{-\infty}^{\infty}(I^d_-\chi_{[0,t]})(s)L(\mathrm{d}s),\qquad t\in\mathbb{R}, \tag{11}$$

where the indicator $\chi_{[a,b]}$ is given by $(a,b\in\mathbb{R})$

$$\chi_{[a,b]}(t)=\begin{cases}1, & \text{if } a\le t<b,\\ -1, & \text{if } b\le t<a,\\ 0, & \text{otherwise.}\end{cases} \tag{12}$$

***Remark 2.2.*** The distribution of $M_d(t)$ is infinitely divisible for all $t \in \mathbb{R}$,

$$\mathrm{E}[M_d(t)^2] = t^{2d+1}\mathrm{E}[L(1)^2], \qquad t \in \mathbb{R}, \quad \text{and}$$

$$\mathrm{E}[\exp\{\mathrm{i}zM_d(t)\}] = \exp\biggl\{\int_{\mathbb{R}}\int_{\mathbb{R}} (\mathrm{e}^{\mathrm{i}zf_t(s)x} - 1 - \mathrm{i}zf_t(s)x)\nu(\mathrm{d}x)\,\mathrm{d}s\biggr\}, \qquad t, z \in \mathbb{R}. \quad (13)$$

## 3. The Lévy Wick exponential and the $S$-transform

One of our aims is to introduce a Hitsuda–Skorokhod integral for convoluted Lévy processes without touching the technicalities of Malliavin calculus and white noise analysis. Our approach is based on the $S$-transform, which uniquely determines a square-integrable random variable by its expectation under an appropriately rich class of probability measures. As a preparation and motivation, we compute the $S$-transform of Itô integrals with respect to the compensated jump measure $\tilde{N}$ in this section. This result then yields a simple definition for anticipative integrals with respect to $\tilde{N}$.

We begin with some definitions.

***Definition 3.1 (Lévy Wick exponential).*** *Let $\mathcal{S}(\mathbb{R}^2)$ denote the Schwartz space of rapidly decreasing smooth functions on $\mathbb{R}^2$. For $\eta \in \Xi$, where*

$$\Xi = \biggl\{\eta \in S(\mathbb{R}^2) : \eta(x,t) > -1, \eta(0,t) = 0, \frac{\mathrm{d}}{\mathrm{d}x}\eta(0,t) = 0, \textit{ for all } t, x \in \mathbb{R}\biggr\},$$

*the* Wiener integral *is defined by*

$$I_1(\eta) = \int_{\mathbb{R}}\int_{\mathbb{R}_0} \eta(x,s)\tilde{N}(\mathrm{d}x,\mathrm{d}s) \quad (14)$$

*and the* Wick exponential *of $I_1(\eta)$ by*

$$\exp^{\diamond}(I_1(\eta)) = \exp\biggl\{\int_{\mathbb{R}}\int_{\mathbb{R}_0} \log(1+\eta(x,t))N(\mathrm{d}x,\mathrm{d}t) - \int_{\mathbb{R}}\int_{\mathbb{R}_0} \eta(x,t)\nu(\mathrm{d}x,\mathrm{d}t)\biggr\}. \quad (15)$$

***Remark 3.2.*** (i) By Theorem 3.1 in Lee and Shih [11],

$$\exp^{\diamond}(I_1(\eta)) = \sum_{n=0}^{\infty} \frac{I_n(\eta^{\otimes n})}{n!}, \quad (16)$$

where $I_n$ denotes the multiple Wiener integral of order $n$ with respect to the compensated Lévy measure. This respresentation justifies the name 'Wick exponential'.

(ii) Since $\exp^{\diamond}(I_1(\eta))$ coincides with the Doléans–Dade exponential of $I_1(\eta)$ at $t = \infty$, it is straightforward that for $\eta, \tilde{\eta} \in \Xi$, we have

$$\mathrm{E}[\exp^{\diamond}(I_1(\eta))] = 1 \quad \text{and} \quad \mathrm{E}[\exp^{\diamond}(I_1(\eta)) \cdot \exp^{\diamond}(I_1(\tilde{\eta}))] = \exp\{(\eta, \tilde{\eta})_{L^2(\nu\times\lambda)}\},$$

where $\lambda$ denotes the Lebesgue measure.

We can now define the $S$-transform.

***Definition 3.3 (S-transform).*** *For $X \in L^2(\Omega, P)$, the S-transform $SX$ of $X$ is an integral transform defined on the set $\Xi$ by*

$$(SX)(\eta) = \mathrm{E}^{Q_\eta}[X], \tag{17}$$

*where*

$$\mathrm{d}Q_\eta = \exp^\diamond(I_1(\eta))\,\mathrm{d}P.$$

Various definitions of the $S$-transform can be found in the literature, which differ according to the chosen subset of deterministic integrands. Our choice of $\Xi$ is particularly convenient because of the smoothness of its members. Moreover, it is a sufficiently rich set, as demonstrated by the following theorem. It states that every square-integrable random variable is uniquely determined by its $S$-transform.

**Proposition 3.4.** *The $S$-transform is injective, that is, if $S(X)(\eta) = S(Y)(\eta)$ for all $\eta \in \Xi$, then $X = Y$.*

**Proof.** The assertion is proved in Løkka and Proske [12], Theorem 5.3, by reformulating a more general result from Albeverio et al. [1], Theorem 5. □

We shall now calculate the $S$-transform of an Itô integral with respect to the compensated jump measure $\tilde{N}$. To this end, let $T > 0$ and $X : \mathbb{R}_0 \times [0, T] \times \Omega \to \mathbb{R}$ be a predictable random field (with respect to the filtration $\mathcal{F}_t$ generated by the Lévy process $L(s), 0 \le s \le t$) satisfying

$$\mathrm{E}\left[\int_0^T \int_{\mathbb{R}_0} |X(y,t)|^2 \nu(\mathrm{d}y)\,\mathrm{d}t\right] < \infty.$$

The compensated Poisson integral $\int_0^T \int_{\mathbb{R}_0} X(y,t)\tilde{N}(\mathrm{d}y, \mathrm{d}t)$ then exists in $L^2(\Omega, P)$.

The following theorem characterizes this integral in terms of the $S$-transform. The result was derived by Løkka and Proske [12], Corollary 7.4, by lengthy calculations involving multiple Wiener integrals. Here, we provide a short proof which only makes use of classical tools such as the Girsanov theorem.

**Theorem 3.5.** *Let $X$ denote a predictable random field satisfying the above integrability condition. $\int_0^T \int_{\mathbb{R}_0} X(y,t)\tilde{N}(\mathrm{d}y, \mathrm{d}t)$ is then the unique square-integrable random variable with $S$-transform given by*

$$\int_0^T \int_{\mathbb{R}_0} S(X(y,t))(\eta)\eta(y,t)\nu(\mathrm{d}y)\,\mathrm{d}t, \qquad \eta \in \Xi. \tag{18}$$

**Proof.** Applying Girsanov's theorem for random measures (Jacod and Shiryaev [8], Theorem 3.17), we obtain that under the measure $Q_\eta$, the compensator of $N(\mathrm{d}y,\mathrm{d}t)$ is given by $(1+\eta(y,t))\nu(\mathrm{d}y)\,\mathrm{d}t$. Hence,

$$\int_0^T \int_{\mathbb{R}_0} X(y,t)\tilde{N}(\mathrm{d}y,\mathrm{d}t) - \int_0^T \int_{\mathbb{R}_0} X(y,t)\eta(y,t)\nu(\mathrm{d}y)\,\mathrm{d}t \tag{19}$$

is a $Q_\eta$-local martingale. In particular, if $0=\tau_1\le\cdots\le\tau_N<\infty$ is a localizing sequence of stopping times with $\lim_{N\to\infty}\tau_N=\infty$ a.s., then

$$\begin{aligned}
&\lim_{N\to\infty} \mathrm{E}^{Q_\eta}\left[\int_0^{T\wedge\tau_N}\int_{\mathbb{R}_0} X(t,y)\tilde{N}(\mathrm{d}y,\mathrm{d}t)\right]\\
&\quad=\lim_{N\to\infty} \mathrm{E}^{Q_\eta}\left[\int_0^{T\wedge\tau_N}\int_{\mathbb{R}_0} X(t,y)\eta(y,t)\nu(\mathrm{d}y)\,\mathrm{d}t\right]\\
&\quad=\mathrm{E}^{Q_\eta}\left[\int_0^{T}\int_{\mathbb{R}_0} X(t,y)\eta(y,t)\nu(\mathrm{d}y)\,\mathrm{d}t\right]
\end{aligned}$$

by a straightforward application of the dominated convergence theorem.

To treat the limit in the first line, note that

$$\mathrm{E}^{Q_\eta}\left[\int_0^{T\wedge\tau_N}\int_{\mathbb{R}_0} X(y,t)\tilde{N}(\mathrm{d}y,\mathrm{d}t)\right]=\mathrm{E}^{P}\left[\exp^\diamond(I_1(\eta))\int_0^{T\wedge\tau_N}\int_{\mathbb{R}_0} X(y,t)\tilde{N}(\mathrm{d}y,\mathrm{d}t)\right].$$

The integrand on the right-hand side is dominated by

$$\exp^\diamond(I_1(\eta))\sup_{0\le u\le T}\left|\int_0^{u}\int_{\mathbb{R}_0} X(y,t)\tilde{N}(\mathrm{d}y,\mathrm{d}t)\right|,$$

which is $P$-integrable by Hölder's inequality, Doob's inequality and the assumed integrability of the random field. Thus,

$$\mathrm{E}^{Q_\eta}\left[\int_0^{T}\int_{\mathbb{R}_0} X(y,t)\tilde{N}(\mathrm{d}y,\mathrm{d}t)\right]=\mathrm{E}^{Q_\eta}\left[\int_0^{T}\int_{\mathbb{R}_0} X(y,t)\eta(y,t)\nu(\mathrm{d}y)\,\mathrm{d}t\right]$$

and the assertion follows by applying Fubini's theorem.

Note that the last identity shows that the $Q_\eta$-local martingale (19) is indeed a $Q_\eta$-martingale. □

***Example 3.6.*** By separately applying the previous theorem to both sides of the two-sided Lévy process, we derive

$$S(M(t))(\eta)=\int_{-\infty}^{t}\int_{\mathbb{R}_0} f(t,s)y\eta(y,s)\nu(\mathrm{d}y)\,\mathrm{d}s$$

since

$$M(t) = \int_{-\infty}^{t} f(t,s)L(\mathrm{d}s) = \int_{-\infty}^{t}\int_{\mathbb{R}_0} f(t,s)y\tilde{N}(\mathrm{d}y,\mathrm{d}s).$$

The $S$-transform characterization in the previous theorem gives rise to a straightforward extension to anticipative random fields.

***Definition 3.7.*** *Suppose $X$ is a random field.*

(i) *The* Hitsuda–Skorokhod integral *of $X$ with respect to the compensated jump measure $\tilde{N}$ is said to exist in $L^2(\Omega)$ if there is a random variable $\Phi \in L^2(\Omega)$ such that for all $\eta \in \Xi$,*

$$S\Phi(\eta) = \int_0^T \int_{\mathbb{R}_0} S(X(y,t))(\eta)\eta(y,t)\nu(\mathrm{d}y)\,\mathrm{d}t.$$

*It is denoted by $\Phi = \int_0^T \int_{\mathbb{R}_0} X(y,t)\tilde{N}^{\diamond}(\mathrm{d}y,\mathrm{d}t)$.*

(ii) *The* Hitsuda–Skorokhod integral *of $X$ with respect to the jump measure $N$ is defined as*

$$\int_0^T \int_{\mathbb{R}_0} X(y,t)N^{\diamond}(\mathrm{d}y,\mathrm{d}t) := \int_0^T \int_{\mathbb{R}_0} X(y,t)\tilde{N}^{\diamond}(\mathrm{d}y,\mathrm{d}t) + \int_0^T \int_{\mathbb{R}_0} X(y,t)\nu(\mathrm{d}y)\,\mathrm{d}t$$

*if both integrals on the right-hand side exist in $L^2(\Omega)$.*

***Remark 3.8.*** From the previous definition, we get immediately that

$$S\Big(\int_0^T \int_{\mathbb{R}_0} X(y,t)N^{\diamond}(\mathrm{d}y,\mathrm{d}t)\Big)(\eta) = \int_0^T \int_{\mathbb{R}_0} S(X(y,t))(\eta)(1+\eta(y,t))\nu(\mathrm{d}y)\,\mathrm{d}t.$$

Clearly, if the integrand is predictable, this Skorokhod integral reduces to the ordinary stochastic integral for random measures and the diamond can be omitted in this case.

***Remark 3.9.*** Theorem 3.5 implies that

$$S(\tilde{N}(A,[0,t]))(\eta) = \int_0^t \int_A \eta(y,s)\nu(\mathrm{d}y)\,\mathrm{d}s.$$

Hence, we can write, in a suggestive notation,

$$S\Big(\int_0^T \int_{\mathbb{R}_0} X(y,t)\tilde{N}^{\diamond}(\mathrm{d}y,\mathrm{d}t)\Big)(\eta) = \int_0^T \int_{\mathbb{R}_0} S(X(y,t))(\eta)S(\tilde{N}(\mathrm{d}y,\mathrm{d}t))(\eta).$$

In view of Example 3.6, Theorem 3.5 can be specialized to integrals with respect to the Lévy process $L$ as follows.

**Corollary 3.10.** *Let $0 \le a \le b$ and $X : [a,b] \times \Omega \to \mathbb{R}$ be a predictable process such that* $\mathrm{E}[\int_a^b |X(t)|^2 \,\mathrm{d}t] < \infty$. $\int_a^b X(s) L(\mathrm{d}s)$ *is then the unique square-integrable random variable with S-transform given by*

$$\int_a^b \int_{\mathbb{R}_0} S(X(t))(\eta) \frac{\mathrm{d}}{\mathrm{d}t} S(L(t))(\eta) \,\mathrm{d}t, \qquad \eta \in \Xi.$$

We close this section with a remark concerning the relationship between the Skorokhod integral with respect to $\tilde{N}$ and ordinary integration.

***Remark 3.11.*** Suppose $u(y,t)$ is a simple random field of the form

$$u(y,t) = F \mathbf{1}_{A \times (a,b]}(y,t), \qquad a < b \in \mathbb{R}, F \in L^2(\Omega),$$

where $A \subset \mathbb{R}_0$ satisfies $\nu(A) < \infty$. Then, in the sense of ordinary integration,

$$\int_{\mathbb{R}} \int_{\mathbb{R}_0} u(y,t) \tilde{N}(\mathrm{d}y, \mathrm{d}t) = F \tilde{N}(A, (a,b]).$$

We want to relate this expression to Skorokhod integration with respect to $\tilde{N}$. We shall suppose that $F = \exp^{\diamond}(I(f))$ for some $f \in L^2(\nu \times \lambda)$, where, in generalization of Definition 3.1,

$$\exp^{\diamond}(I(f)) := \exp\left\{ \int_{\mathbb{R}} \int_{\mathbb{R}_0} f(y,t) \tilde{N}(\mathrm{d}y, \mathrm{d}t) \right\} \prod_{s \in \mathbb{R}} (1 + f(\Delta L(s), s)) \mathrm{e}^{-f(\Delta L(s), s)}.$$

A direct calculation then shows that for $f \in L^2(\nu \times \lambda)$, $\eta \in \Xi$,

$$\exp^{\diamond}(I(f)) \exp^{\diamond}(I(\eta)) = \exp^{\diamond}(I(f + \eta + f\eta)) \exp\left\{ \int_{\mathbb{R}} \int_{\mathbb{R}_0} f(y,t) \eta(y,t) \nu(\mathrm{d}y) \,\mathrm{d}t \right\}.$$

Consequently, by a slight extension of Remark 3.9,

$$\begin{aligned}
&S\left( \int_{\mathbb{R}} \int_{\mathbb{R}_0} u(y,t) \tilde{N}(\mathrm{d}y, \mathrm{d}t) \right)(\eta) \\
&\quad = \exp\left\{ \int_{\mathbb{R}} \int_{\mathbb{R}_0} f(y,t) \eta(y,t) \nu(\mathrm{d}y) \,\mathrm{d}t \right\} \mathrm{E}[\exp^{\diamond}(I(f + \eta + f\eta)) \tilde{N}(A, (a,b])] \\
&\quad = (S \exp^{\diamond}(I(f)))(\eta) \int_a^b \int_A [\eta(y,s) + f(y,s) + \eta(y,s) f(y,s)] \nu(\mathrm{d}y) \,\mathrm{d}s \\
&\quad = S\left( \int_{\mathbb{R}} \int_{\mathbb{R}_0} u(y,t) f(y,t) \nu(\mathrm{d}y) \,\mathrm{d}t \right)(\eta) \\
&\qquad + S\left( \int_{\mathbb{R}} \int_{\mathbb{R}_0} [u(y,t) f(y,t) + u(y,t)] \tilde{N}^{\diamond}(\mathrm{d}y, \mathrm{d}t) \right)(\eta).
\end{aligned}$$

Let us now define the *Malliavin derivative* of a Wick exponential by

$$D_{y,s}\exp^{\diamond}(I(f)) = f(y,s)\exp^{\diamond}(I(f)),$$

which can be extended to a linear closed operator acting on a larger class of random variables (see, e.g., Nualart and Vives [16]). We then arrive at the formula

$$\begin{aligned}\int_{\mathbb{R}}\int_{\mathbb{R}_0} u(y,t)\tilde{N}(\mathrm{d}y,\mathrm{d}t) = &\int_{\mathbb{R}}\int_{\mathbb{R}_0} D_{y,t}u(y,t)\nu(\mathrm{d}y)\,\mathrm{d}s\\ &+\int_{\mathbb{R}}\int_{\mathbb{R}_0}[u(y,t)+D_{y,t}u(y,t)]\tilde{N}^{\diamond}(\mathrm{d}y,\mathrm{d}s).\end{aligned}$$

We conjecture that this formula can be extended by approximation to a larger class of random fields.

## 4. A Skorokhod integral for convoluted Lévy processes

In this section, we define the Skorokhod integral for convoluted Lévy processes and state some basic properties. The definition is strongly motivated by Corollary 3.10 above.

***Definition 4.1.*** *Suppose that the mapping*

$$t \mapsto S(M(t))(\eta)$$

*is differentiable for every $\eta \in \Xi$. Suppose $B \subset \mathbb{R}$ is a Borel set and $X : B \times \Omega \to \mathbb{R}$ is a measurable stochastic process such that $X(t)$ is square-integrable for each $t \in B$. $X$ is said then to have a* Hitsuda-Skorokhod integral *with respect to $M$ if*

$$S(X(\cdot))(\eta)\frac{\mathrm{d}}{\mathrm{d}t}S(M(\cdot))(\eta) \in L^1(B) \qquad \textit{for any } \eta \in \Xi$$

*and there is a $\Phi \in L^2(\Omega)$ such that for all $\eta \in \Xi$,*

$$S(\Phi)(\eta) = \int_B S(X(t))(\eta)\frac{\mathrm{d}}{\mathrm{d}t}S(M(t))(\eta)\,\mathrm{d}t.$$

In that case, $\Phi$ is uniquely determined by the injectivity of the $S$-transform and we write

$$\Phi = \int_B X(t)M^{\diamond}(\mathrm{d}t).$$

***Remark 4.2.*** (i) Lemma 5.1(ii) below provides some sufficient conditions for the differentiability of the mapping $t \mapsto S(M(t))(\eta)$ in terms of the convolution kernel.

(ii) The definition of the Skorokhod integral does not require conditions such as predictability or progressive measurability. Hence, it also generalizes the Itô integral with respect to the underlying Lévy process to anticipative integrands.

(iii) Since the Lévy process itself is stochastically continuous, the $S$-transform cannot distinguish between $L(t)$ and $L(t-)$ for fixed $t$. Consequently, we obtain, for example,

$$\int_0^t L(s)L^\diamond(\mathrm{d}s) = \int_0^t L(s-)L^\diamond(\mathrm{d}s) = \int_0^t L(s-)L(\mathrm{d}s),$$

where the last integral is the classical Itô integral.

The following properties of the Skorokhod integral are an obvious consequence of the definition.

**Proposition 4.3.** (i) *For all* $a < b \in \mathbb{R}$, $M(b) - M(a) = \int_a^b M^\diamond(\mathrm{d}t)$.

(ii) *Let* $X : B \times \Omega \to L^2(\Omega)$ *be Skorokhod integrable. Then*

$$\int_B X(t)M^\diamond(\mathrm{d}t) = \int_{\mathbb{R}} 1_B(t)X(t)M^\diamond(\mathrm{d}t),$$

*where* $1_B$ *denotes the indicator function of the set* $B$.

(iii) *Let* $X : B \times \Omega \to L^2(\Omega)$ *be Skorokhod integrable. Then* $\mathrm{E}[\int_B X(t)M^\diamond(\mathrm{d}t)] = 0$.

We note that (iii) holds since the expectation coincides with the $S$-transform at $\eta = 0$. The zero expectation property makes the integral a promising candidate for modeling an additive noise.

***Example 4.4.*** As an example, we show how to calculate $\int_0^T M(t)M^\diamond(\mathrm{d}t)$. In the following manipulations, $\tilde{N}^\eta$ denotes the compensated jump measure under the probability measure $Q_\eta = \exp^\diamond(I_1(\eta))\,\mathrm{d}P$. In particular, it follows from Girsanov's theorem, as in the proof of Theorem 3.5, that

$$M(T) = \int_{-\infty}^T \int_{\mathbb{R}_0} f(T,s)y\tilde{N}^\eta(\mathrm{d}y,\mathrm{d}s) + \int_{-\infty}^T \int_{\mathbb{R}_0} f(T,s)y\eta(y,s)\nu(\mathrm{d}y)\,\mathrm{d}s.$$

By this identity, integration by parts and Example 3.6, we obtain

$$\begin{aligned}
&S\left(2\int_0^T M(t)M^\diamond(\mathrm{d}t)\right)(\eta) \\
&\quad = 2\int_0^T S(M(t))(\eta)\frac{\mathrm{d}}{\mathrm{d}t}S(M(t))(\eta)\,\mathrm{d}t \\
&\quad = (S(M(T))(\eta))^2 = \left(\int_{-\infty}^T \int_{\mathbb{R}_0} f(T,s)y\eta(y,s)\nu(\mathrm{d}y)\,\mathrm{d}s\right)^2
\end{aligned}$$

$$= \mathrm{E}^{Q_\eta}\left[\left(\int_{-\infty}^{T}\int_{\mathbb{R}_0} f(T,s)y\tilde{N}^\eta(\mathrm{d}y,\mathrm{d}s) + \int_{-\infty}^{T}\int_{\mathbb{R}_0} f(T,s)y\eta(y,s)\nu(\mathrm{d}y)\,\mathrm{d}s\right)^2\right]$$
$$- \mathrm{E}^{Q_\eta}\left[\left(\int_{-\infty}^{T}\int_{\mathbb{R}_0} f(T,s)y\tilde{N}^\eta(\mathrm{d}y,\mathrm{d}s)\right)^2\right]$$
$$= S(M(T)^2)(\eta) - \int_{-\infty}^{T}\int_{\mathbb{R}_0} f(T,s)^2y^2(1+\eta(y,s))\nu(\mathrm{d}y)\,\mathrm{d}s.$$

Here, we have used the fact that $\int_{-\infty}^{T}\int_{\mathbb{R}_0} f(T,s)y\tilde{N}^\eta(\mathrm{d}y,\mathrm{d}s)$ has zero expectation and variance $\int_{-\infty}^{T}\int_{\mathbb{R}_0} f(T,s)^2y^2(1+\eta(y,s))\nu(\mathrm{d}y)\,\mathrm{d}s$ since the compensator of $N$ under $Q_\eta$ is given by $(1+\eta(y,s))\nu(\mathrm{d}y,\mathrm{d}s)$.

Hence, from Remark 3.8, we derive the identity

$$2\int_0^T M(t)M^\diamond(\mathrm{d}t) = M(T)^2 - \int_{-\infty}^{T}\int_{\mathbb{R}_0} f(T,s)^2y^2N(\mathrm{d}y,\mathrm{d}s)$$
$$= M(T)^2 - \sum_{-\infty<s\le T} f(T,s)^2(\Delta L(s))^2,$$

provided both members on the right-hand side exist in $L^2(\Omega)$.

***Remark 4.5.*** Applying the same techniques as in Example 4.4, one can easily obtain, for $a \le b$,

$$M(a)\int_a^b 1M^\diamond(\mathrm{d}t) = \int_a^b M(a)M^\diamond(\mathrm{d}t) + \int_0^a\int_{\mathbb{R}_0} f(a,s)(f(b,s)-f(a,s))y^2N(\mathrm{d}y,\mathrm{d}s).$$

Hence, ordinary multiplication with a random variable, which is measurable with respect to the information up to the lower integration bound, cannot, in general, be introduced under the integral sign if the kernel depends on the past.

## 5. Itô's formula

In this section, we will derive an Itô formula for convoluted Lévy processes. The proof is based on a calculation of the time derivative of $S(G(M(t)))(\eta)$. It may be seen as a generalization of the calculations in Example 4.4. This technique of proof is in the spirit of Kubo [9], Bender [3] and Lee and Shih [10], where this approach was applied to obtain Itô formulas for generalized functionals of a Brownian motion a fractional Brownian motion, and a Lévy process with Brownian component, respectively.

During the derivation of the Itô formula, we have to interchange differentiation and integration several times. Under the following (rather strong) conditions on the convolution kernel, these manipulations are easily justified. However, the Itô formulas below may also be viewed as generic results which hold for more general kernels (with the technicalities to be checked on a case-by-case basis).

We make the following assumptions:

(H1) there are constants $a \leq 0 < b$ such that $\operatorname{supp}(f) \subset [a,b]^2$;
(H2) $f$ is continuous and bounded on $[a,b]^2 \setminus \{(t,s); t=s\}$;
(H3) $\lim_{s\uparrow t} f(t,s) = f(t,t)$ and the mapping $t \mapsto f(t,t)$ is continuous;
(H4) $f$ is continuously differentiable on $(a,b)^2 \setminus \{(t,s); t=s\}$ with bounded derivative.

**Lemma 5.1.** *Under* (H1)–(H4), *we have the following:*

(i) *For* $a \leq t \leq b$,

$$M(t) = f(t,t)L(t) - f(t,a)L(a) - \int_a^t L(s)\frac{\mathrm{d}}{\mathrm{d}s}f(t,s)\,\mathrm{d}s. \tag{20}$$

*In particular,* $M(t)$ *has a modification which is RCLL and stochastically continuous. Moreover,*

$$\Delta M(t) = f(t,t)\Delta L(t).$$

*Hence,* $M$ *is continuous on* $[a,b]$ *if and only if* $f(t,t)=0$ *for all* $a \leq t \leq b$.

(ii) *The mapping* $[a,b] \to \mathbb{R}, t \mapsto (SM(t))(\eta)$ *is continuously differentiable for all* $\eta \in \Xi$ *and*

$$\frac{\mathrm{d}}{\mathrm{d}t}(SM(t))(\eta) = \int_{-\infty}^t \int_{\mathbb{R}_0} \frac{\mathrm{d}}{\mathrm{d}t} f(t,s) y\eta(y,s)\nu(\mathrm{d}y)\,\mathrm{d}s + f(t,t)\int_{\mathbb{R}_0} y\eta(y,t)\nu(\mathrm{d}y). \tag{21}$$

**Proof.** (i) Fix a modification of $L$ which is right-continuous with left limits (RCLL). Formula (20) follows from the definition of $M$ and integration by parts, which is justified by (H3)–(H4). The second and third terms on the right-hand side are continuous in $t$ by (H2) and the boundedness of $\frac{\mathrm{d}}{\mathrm{d}s}f(t,s)$, respectively. The first term is stochastically continuous and RCLL since $L$ has these properties and $t \mapsto f(t,t)$ is continuous. The other assertions in (i) are obvious consequences.

(ii) can easily be obtained by differentiating the expression in Example 3.6. □

***Example 5.2.*** The following prominent examples satisfy conditions (H1)–(H4):

1. *one-sided shot noise processes* defined by the kernel

$$f(t,s) = \begin{cases} k(t-s), & 0 \leq s \leq t \leq T^*, \\ 0, & \text{otherwise,} \end{cases}$$

   for constants $T^* > 0$ and $k$;
2. *one-sided Ornstein–Uhlenbeck type processes* defined by the kernel

$$f(t,s) = \begin{cases} \mathrm{e}^{-k(t-s)}, & 0 \leq s \leq t \leq T^*, \\ 0, & \text{otherwise,} \end{cases}$$

for constants $T^* > 0$ and $k \geq 0$.

From the previous lemma, one directly obtains that the shot noise processes have continuous paths, while the Ornstein–Uhlenbeck-type processes exhibit jumps.

To state the Itô formula precisely, we finally recall that the *Wiener algebra* is defined as

$$A(\mathbb{R}) := \{G \in L^1(\mathbb{R}); \mathcal{F}G \in L^1(\mathbb{R})\},$$

where $\mathcal{F}$ denotes the Fourier transform. Note that the space of rapidly decreasing smooth functions is included in the Wiener algebra.

The first version of Itô's formula requires that the underlying Lévy process is a finite variation process.

**Theorem 5.3 (Itô formula I).** *Let* (H1)–(H4) *hold,* $0 < T \leq b$ *and*

$$\int_{\mathbb{R}_0} |x|\nu(\mathrm{d}x) < \infty.$$

*Furthermore, assume that* $G \in C^1(\mathbb{R})$ *with* $G, G' \in A(\mathbb{R})$ *bounded. Then*

$$\int_0^T \left( \int_{-\infty}^t \int_{\mathbb{R}_0} G'(M(t) + xf(t,s)) x \frac{\mathrm{d}}{\mathrm{d}t} f(t,s) N^\diamond(\mathrm{d}x, \mathrm{d}s) \right) \mathrm{d}t$$

*exists in* $L^2(\Omega)$ *and*

$$\begin{aligned} G(M(T)) = G(0) &- \left( \int_{\mathbb{R}_0} x\nu(\mathrm{d}x) \right) \int_0^T G'(M(t)) \left( f(t,t) + \int_{-\infty}^t \frac{\mathrm{d}}{\mathrm{d}t} f(t,s)\, \mathrm{d}s \right) \mathrm{d}t \\ &+ \sum_{0 \leq t \leq T} G(M(t)) - G(M(t-)) \\ &+ \int_0^T \left( \int_{-\infty}^t \int_{\mathbb{R}_0} G'(M(t-) + xf(t,s)) x \frac{\mathrm{d}}{\mathrm{d}t} f(t,s) N^\diamond(\mathrm{d}x, \mathrm{d}s) \right) \mathrm{d}t. \end{aligned}$$

In the general case, the Itô formula reads as follows. Indeed, the previous formula can be derived from the general one by rearranging some terms.

**Theorem 5.4 (Itô formula II).** *Let* (H1)–(H4) *hold and* $0 < T \leq b$. *Furthermore, assume that* $G \in C^1(\mathbb{R})$ *with* $G, G' \in A(\mathbb{R})$. *Then*

$$\begin{aligned} &G(M(T)) \\ &\quad = G(0) + \int_0^T G'(M(t-)) M^\diamond(\mathrm{d}t) \end{aligned}$$

$$+\sum_{0\le t\le T} G(M(t)) - G(M(t-)) - G'(M(t-))\Delta M(t)$$

$$+\int_0^T \left(\int_{-\infty}^t \int_{\mathbb{R}_0} (G'(M(t-)+xf(t,s)) - G'(M(t-)))x\frac{\mathrm{d}}{\mathrm{d}t}f(t,s)N^{\diamond}(\mathrm{d}x,\mathrm{d}s)\right)\mathrm{d}t,$$

*provided all terms exist in* $L^2(\Omega)$.

We would like to emphasize that the Skorokhod integrals with respect to $N$ in the above versions of Itô's formula do not, in general, reduce to ordinary integrals for the following reason. The time variable of the Skorokhod integral is $s$, but the integrand depends on $M$ through the value $M(t-)$, where $t > s$. Therefore, the integrand is typically not predictable as a process in the variable $s$.

The above versions of Itô's formula (but not their exact assumptions) reduce to well-known formulas for Lévy processes with the choice $f(t,s) = \chi_{(0,t]}(s)$ as, in this case, the last Skorokhod integral with respect to $N$ vanishes. We recall that $M$ has independent increments if and only if $\frac{\mathrm{d}}{\mathrm{d}t}f(t,s) = 0$ for all $t$. Hence, the contributions from discontinuities and memory effects are nicely separated in the above Itô formulas. Finally, note that the formula for $M(t)^2$ from Example 4.4 can be recovered by formally applying the Itô formula II with $G(y) = y^2$.

***Remark 5.5.*** Itô formula II has the drawback that the conditions do not guarantee that all members of the identity exist in $L^2(\Omega)$. However, the manipulations below can be recast in a white noise framework, as developed in [17], in a way that all members exist as generalized random variables.

The remainder of this section is devoted to the proof of the Itô formulas. As a general strategy, we wish to show that both sides of the asserted identities have the same $S$-transform. Indeed, the following calculations show how to identify the right-hand side constructively. We first write

$$S(G(M(T)))(\eta) = G(0) + \int_0^T \frac{\mathrm{d}}{\mathrm{d}t} S(G(M(t)))(\eta)\,\mathrm{d}t$$

and then calculate $\frac{\mathrm{d}}{\mathrm{d}t}S(G(M(t)))$ explicitly. To achieve this, we apply the inverse Fourier theorem and obtain, for $G \in A(\mathbb{R})$,

$$S(G(M(t)))(\eta) = \mathrm{E}^{Q_\eta}[G[M(t)]] = \frac{1}{\sqrt{2\pi}} \int_{\mathbb{R}} \mathcal{F}G(u)\mathrm{E}^{Q_\eta}[\mathrm{e}^{\mathrm{i}uM(t)}]\,\mathrm{d}u. \tag{22}$$

To differentiate this expression, we calculate the characteristic function of $M$ under $Q_\eta$.

**Proposition 5.6.** *Let* $M = \{M(t)\}_{t\in\mathbb{R}}$ *be a convoluted Lévy process as defined in* (6), *with kernel function* $f$. *Then*

$S(\mathrm{e}^{\mathrm{i}uM(t)})(\eta)$

$$= \mathrm{E}^{Q_\eta}[\mathrm{e}^{\mathrm{i}uM(t)}]$$
$$= \exp\left\{\mathrm{i}uS(M(t))(\eta) + \int_{-\infty}^{t}\int_{\mathbb{R}_0} [(\mathrm{e}^{\mathrm{i}uxf(t,s)} - 1 - \mathrm{i}uxf(t,s))(1+\eta(x,s))]\nu(\mathrm{d}x)\,\mathrm{d}s\right\}.$$

**Proof.** It follows from the proof of Theorem 3.5 that

$$L^Q(t) := L(t) - \int_0^t \int_{\mathbb{R}_0} x\eta(x,s)\nu(\mathrm{d}x)\,\mathrm{d}s$$

is a $Q_\eta$-martingale with zero mean. Applying Girsanov's theorem for semimartingales (Jacod and Shiryayev [8], Theorem 3.7) yields that $L^Q$ has semimartingale characteristics $(\gamma_s^Q, 0, \nu_s^Q)$, where $\gamma_s^Q = -\int_{|x|>1} x(1+\eta(x,s))\nu(\mathrm{d}x)$ and $\nu_s^Q(\mathrm{d}x) = (1+\eta(x,s))\nu(\mathrm{d}x)$. Hence,

$$S(\exp\{\mathrm{i}uL^Q(t)\})(\eta) = \exp\left\{\int_0^t \int_{\mathbb{R}_0} [\mathrm{e}^{\mathrm{i}ux} - 1 - \mathrm{i}ux][1+\eta(x,s)]\nu(\mathrm{d}x)\,\mathrm{d}s\right\}.$$

Finally,

$$\begin{aligned}
&S(\exp\{\mathrm{i}uM(t)\})(\eta)\\
&\quad = \mathrm{E}^{Q_\eta}\left[\exp\left\{\mathrm{i}u\int_{-\infty}^t f(t,s)L(\mathrm{d}s)\right\}\right]\\
&\quad = \mathrm{E}^{Q_\eta}\left[\exp\left\{\mathrm{i}u\int_{-\infty}^t f(t,s)L^Q(\mathrm{d}s) + \mathrm{i}u\int_{-\infty}^t f(t,s)\int_{\mathbb{R}_0} x\eta(x,s)\nu(\mathrm{d}x)\,\mathrm{d}s\right\}\right]\\
&\quad = \exp\left\{\int_{-\infty}^t \int_{\mathbb{R}_0} [\mathrm{e}^{\mathrm{i}uxf(t,s)} - 1 - \mathrm{i}uxf(t,s)][1+\eta(x,s)]\nu(\mathrm{d}x)\,\mathrm{d}s\right\}\\
&\qquad \times \exp\left\{\int_{-\infty}^t \int_{\mathbb{R}_0} \mathrm{i}uxf(t,s)\eta(x,s)\nu(\mathrm{d}x)\,\mathrm{d}s\right\}.
\end{aligned}$$

Taking the $S$-transform of $M$ into account, which was calculated in Example 3.6, the assertion follows. □

By introducing the derivative under the integral sign, we get

$$\begin{aligned}
&\frac{\mathrm{d}}{\mathrm{d}t}\mathrm{E}^{Q_\eta}[\mathrm{e}^{\mathrm{i}uM(t)}]\\
&\quad = \mathrm{E}^{Q_\eta}[\mathrm{e}^{\mathrm{i}uM(t)}]\int_{\mathbb{R}_0} [(\mathrm{e}^{\mathrm{i}uxf(t,t)} - 1 - \mathrm{i}uxf(t,t))(1+\eta(x,t))]\nu(\mathrm{d}x)\\
&\qquad + \mathrm{E}^{Q_\eta}[\mathrm{e}^{\mathrm{i}uM(t)}]\int_{-\infty}^t \int_{\mathbb{R}_0} \left[\mathrm{i}ux\frac{\mathrm{d}}{\mathrm{d}t}f(t,s)(\mathrm{e}^{\mathrm{i}uxf(t,s)} - 1)(1+\eta(x,t))\right]\nu(\mathrm{d}x)\,\mathrm{d}s\\
&\qquad + \mathrm{E}^{Q_\eta}[\mathrm{e}^{\mathrm{i}uM(t)}]\mathrm{i}u\frac{\mathrm{d}}{\mathrm{d}t}S(M(t))(\eta).
\end{aligned} \tag{23}$$

Combining (22) with (23) and again interchanging differentiation and integration (which can be justified under (H1)–(H4) since $G, G' \in A(\mathbb{R})$), we obtain

$$\begin{aligned}
&\frac{\mathrm{d}}{\mathrm{d}t} S(G(M(t)))(\eta) \\
&\quad = \frac{1}{\sqrt{2\pi}} \int_{\mathbb{R}} (\mathcal{F}G)(u) \mathrm{E}^{Q_\eta}[\mathrm{e}^{\mathrm{i}uM(t)}] \int_{\mathbb{R}_0} [(\mathrm{e}^{\mathrm{i}uxf(t,t)} - 1 - \mathrm{i}uxf(t,t))(1+\eta(x,t))]\nu(\mathrm{d}x)\,\mathrm{d}u \\
&\qquad + \frac{1}{\sqrt{2\pi}} \int_{\mathbb{R}} (\mathcal{F}G)(u) \mathrm{E}^{Q_\eta}[\mathrm{e}^{\mathrm{i}uM(t)}] \\
&\qquad\qquad \times \int_{-\infty}^{t} \int_{\mathbb{R}_0} \left[\mathrm{i}ux \frac{\mathrm{d}}{\mathrm{d}t} f(t,s)(\mathrm{e}^{\mathrm{i}uxf(t,s)} - 1)(1+\eta(x,t))\right] \nu(\mathrm{d}x)\,\mathrm{d}s\,\mathrm{d}u \\
&\qquad + \frac{1}{\sqrt{2\pi}} \int_{\mathbb{R}} (\mathcal{F}G)(u) \mathrm{E}^{Q_\eta}[\mathrm{e}^{\mathrm{i}uM(t)}] \mathrm{i}u \frac{\mathrm{d}}{\mathrm{d}t} S(M(t))(\eta)\,\mathrm{d}u \\
&\quad =: (I) + (II) + (III).
\end{aligned}$$

Standard manipulations of the Fourier transform, together with (22), now yield

$$\begin{aligned}
(I) &= \frac{1}{\sqrt{2\pi}} \int_{\mathbb{R}} \int_{\mathbb{R}_0} [(\mathcal{F}G(\cdot + xf(t,t))(u) - (\mathcal{F}G)(u) - xf(t,t)(\mathcal{F}G')(u)] \\
&\qquad\qquad \times \mathrm{E}^{Q_\eta}[\mathrm{e}^{\mathrm{i}uM(t)}](1+\eta(x,t))\nu(\mathrm{d}x)\,\mathrm{d}u \\
&= \int_{\mathbb{R}_0} S(G(M(t-) + xf(t,t)) - G(M(t-)) - xf(t,t)G'(M(t-)))(\eta) \\
&\qquad \times (1+\eta(x,t))\nu(\mathrm{d}x).
\end{aligned}$$

The second term can be treated analogously and thus,

$$\begin{aligned}
(II) = \int_{-\infty}^{t} \int_{\mathbb{R}_0} x \frac{\mathrm{d}}{\mathrm{d}t} f(t,s) S(G'(M(t-) + xf(t,s)) - G'(M(t-)))(\eta) \\
\times (1+\eta(x,t))\nu(\mathrm{d}x)\,\mathrm{d}s.
\end{aligned}$$

Finally, $(III) = S(G'(M(t-)))(\eta) \frac{\mathrm{d}}{\mathrm{d}t} S(M(t))(\eta)$.

We now collect terms and integrate $t$ from 0 to $T$, whence

$$\begin{aligned}
&S(G(M(T)))(\eta) - G(0) \\
&\quad = \int_0^T \int_{\mathbb{R}_0} S(G(M(t-) + xf(t,t)) - G(M(t-)) - xf(t,t)G'(M(t-)))(\eta) \\
&\qquad\qquad \times (1+\eta(x,t))\nu(\mathrm{d}x)\,\mathrm{d}t \\
&\qquad + \int_0^T \int_{-\infty}^{t} \int_{\mathbb{R}_0} x \frac{\mathrm{d}}{\mathrm{d}t} f(t,s) S(G'(M(t-) + xf(t,s)) - G'(M(t-)))(\eta)
\end{aligned}$$

$$\times (1+\eta(x,t))\nu(\mathrm{d}x)\,\mathrm{d}s\,\mathrm{d}t$$
$$+\int_0^T S(G'(M(t-)))(\eta)\frac{\mathrm{d}}{\mathrm{d}t}S(M(t))(\eta)\,\mathrm{d}t$$
$$=: (i)+(ii)+(iii). \tag{24}$$

From Remark 3.8, we get

$$\begin{aligned}(i) &= S\bigg(\int_0^T\int_{\mathbb{R}_0} G(M(t-)+xf(t,t))-G(M(t-))-xf(t,t)G'(M(t-))N^\diamond(\mathrm{d}x,\mathrm{d}t)\bigg)(\eta)\\ &= S\Bigg(\sum_{0\le t\le T} G(M(t))-G(M(t-))-G'(M(t-))\Delta M(t)\Bigg)(\eta),\end{aligned}$$

where the second identity holds because the Skorokhod integral is an Itô integral by predictability (and by taking Lemma 5.1(i) into account). Similarly,

$$(ii) = S\bigg(\int_0^T\int_{-\infty}^t\int_{\mathbb{R}_0} x\frac{\mathrm{d}}{\mathrm{d}t}f(t,s)[G'(M(t-)+xf(t,s))-G'(M(t-))]N^\diamond(\mathrm{d}x,\mathrm{d}s)\,\mathrm{d}t\bigg)(\eta).$$

Finally, by the definition of the Skorokhod integral with respect to $M$,

$$(iii) = S\bigg(\int_0^T G'(M(t-))M^\diamond(\mathrm{d}t)\bigg)(\eta).$$

Hence, both sides of Itô formula II have the same $S$-transform, which proves this formula.

To get Itô formula I, we rearrange the terms in (24). By Lemma 5.1(ii),

$$\frac{\mathrm{d}}{\mathrm{d}t}S(M(t))(\eta) = f(t,t)\int_{\mathbb{R}_0} x\eta(x,t)\nu(\mathrm{d}x)+\int_{-\infty}^t\frac{\mathrm{d}}{\mathrm{d}t}f(t,s)\int_{\mathbb{R}_0} x\eta(x,s)\nu(\mathrm{d}x)\,\mathrm{d}s.$$

Thus, by (24) and similar considerations as above,

$$\begin{aligned}&\int_0^T\int_{-\infty}^t\int_{\mathbb{R}_0} x\frac{\mathrm{d}}{\mathrm{d}t}f(t,s)S(G'(M(t-)+xf(t,s)))(\eta)(1+\eta(x,t))\nu(\mathrm{d}x)\,\mathrm{d}s\,\mathrm{d}t\\ &\quad= S\Bigg(G(M(T))-G(0)+\bigg(\int_{\mathbb{R}_0} x\nu(\mathrm{d}x)\bigg)\int_0^T G'(M(t))\bigg(f(t,t)+\int_{-\infty}^t\frac{\mathrm{d}}{\mathrm{d}t}f(t,s)\,\mathrm{d}s\bigg)\mathrm{d}t\\ &\qquad -\sum_{0\le t\le T} G(M(t-)+\Delta M(t))-G(M(t-))\Bigg)(\eta).\end{aligned}$$

The expression under the $S$-transform on the right-hand side clearly belongs to $L^2(\Omega)$ under the assumptions of Itô formula I. Then, by Remark 3.8, the Skorokhod integral

$$\int_0^T \int_{-\infty}^t \int_{\mathbb{R}_0} x \frac{\mathrm{d}}{\mathrm{d}t} f(t,s) G'(M(t-) + xf(t,s)) N^{\diamond}(\mathrm{d}x, \mathrm{d}s)\, \mathrm{d}t$$

exists in $L^2(\Omega)$ and coincides with the expression under the $S$-transform on the right-hand side. This proves Itô formula I.

# 6. Stochastic calculus for fractional Lévy processes

We shall now specialize from a convoluted Lévy process to a fractional one. In Marquardt [14], a Wiener-type integral with respect to a fractional Lévy process is defined for deterministic integrands. Its domain is the space of functions $g$ such that $I_-^d g \in L^2(\mathbb{R})$ and it can be characterized by the property

$$\int_{\mathbb{R}} g(s) M_d(\mathrm{d}s) = \int_{\mathbb{R}} (I_-^d g)(s) L(\mathrm{d}s).$$

The following theorem shows that a similar characterization holds for Skorokhod integrals with respect to fractional Lévy processes. Hence, it also proves, as a by-product that the Wiener-type integral is a special case of the Skorokhod integral.

In preparation, note that

$$S(M_d(t))(\eta) = \int_{\mathbb{R}} \int_{\mathbb{R}_0} I_-^d \chi_{[0,t]}(s) y \eta(s,y) \nu(\mathrm{d}y)\, \mathrm{d}s.$$

Hence, by Fubini's theorem and fractional integration by parts, we obtain the following theorem.

**Theorem 6.1.** *Suppose $M_d$ is a fractional Lévy process with $0 < d < 0.5$. Then, for all $\eta \in \Xi$,*

$$\frac{\mathrm{d}}{\mathrm{d}t} S(M_d(t))(\eta) = \int_{\mathbb{R}_0} (I_+^d \eta)(t,y) y \nu(\mathrm{d}y),$$

*where, by convention, fractional integral operators are applied only to the time variable $t$.*

*Furthermore, suppose that $X \in L^p(\mathbb{R}, L^2(\Omega))$ with $p = (d + 1/2)^{-1}$. Then*

$$\int_{\mathbb{R}} X(t) M_d^{\diamond}(\mathrm{d}t) = \int_{\mathbb{R}} (I_-^d X)(t) L^{\diamond}(\mathrm{d}t)$$

*in the usual sense, that is, if one of the integrals exists, then so does the other and both coincide.*

**Proof.** The proof follows the same lines as that of Theorem 3.4 in Bender [4]. □

Note that only Itô formula II makes sense for fractional Lévy processes. When we formally apply this Itô formula, the following observation is noteworthy. For $d>0$, the process $M_d$ is continuous and has memory, whence

$$\begin{aligned}G(M_d(T)) = G(0) + \int_0^T G'(M_d(t-))M^\diamond(\mathrm{d}t)\\ + \int_0^T \Big( \int_{-\infty}^t \int_{\mathbb{R}_0} \Big( G'\Big( M_d(t-) + \frac{x}{\Gamma(d+1)}((t-s)_+^d - (-s)_+^d) \Big)\\ - G'(M_d(t-)) \Big) \frac{x}{\Gamma(d)} (t-s)_+^{d-1} N^\diamond(\mathrm{d}x, \mathrm{d}s) \Big) \, \mathrm{d}t.\end{aligned}$$

However, the Lévy process $L$ itself comes up as limit of $M_d$ when $d$ tends to 0. As this process has independent increments and jumps, its well-known Itô formula reads

$$\begin{aligned}G(L(T)) = G(0) + \int_0^T G'(L(t-))L(\mathrm{d}t)\\ + \sum_{0\le t\le T} G(L(t)) - G(L(t-)) - G'(L(t-))\Delta L(t).\end{aligned}$$

So, apparently, the Itô formulas do not transform continuously into each other when passing to this limit. This is in sharp contrast to the Gaussian case, in which the Itô formula for Brownian motion is recovered by substituting $H=1/2$ (the Hurst parameter corresponding to $d$ via $d=H-1/2$) into the Itô formula for fractional Brownian motions (see, e.g., [4]).

## Acknowledgements

The paper benefited from the constructive comments of two anonymous referees. In particular, their remarks helped to clarify the conditions required for the proof of Itô's formula.